\documentclass[graybox]{svmult}

\usepackage{type1cm}        
\usepackage{makeidx}         
\usepackage{graphicx}        
\usepackage{multicol}        
\usepackage[bottom]{footmisc}
\usepackage{caption}

\usepackage{newtxtext}       %
\usepackage[varvw]{newtxmath}       

\usepackage{amsmath}
\usepackage{amsfonts}
\usepackage{amsthm}
\usepackage{tikz}
\usepackage{xcolor}
\usepackage{newtxtext,newtxmath}
\usepackage{caption}
\usepackage{subcaption}
\usepackage{hyperref}
\usepackage{todonotes}
\usepackage{soul}

\makeindex             

\begin{document}

\title*{Optimisation--Based Coupling of Finite Element Model and Reduced Order Model for Computational Fluid Dynamics}
\titlerunning{Optimisation--Based Coupling of FEM and ROM for CFD}
\author{Ivan Prusak\orcidID{0009-0008-2882-0576}  and\\ Davide Torlo\orcidID{0000-0001-5106-1629}   and\\ Monica Nonino\orcidID{0000-0002-5503-705X}  and\\ Gianluigi Rozza\orcidID{0000-0002-0810-8812}}
\institute{Ivan Prusak \at  Sissa, Mathematics Area, mathLab, International School for Advanced Studies, via Bonomea 265, 34136 Trieste, Italy, \email{iprusak@sissa.it}
\and Davide Torlo \at  Università di Roma La Sapienza, Dipartimento di Matematica ``Guido Castelnuovo'', piazzale Aldo Moro 5, 00185, Rome, Italy, \email{davide.torlo@uniroma1.it} 
\and Monica Nonino \at University of Vienna, Department of Mathematics, Oskar-Morgenstern-Platz 1, 1090 Vienna, Austria, \email{monica.nonino@univie.ac.at}
\and Gianluigi Rozza \at  Sissa, Mathematics Area, mathLab, International School for Advanced Studies, via Bonomea 265, 34136 Trieste, Italy, \email{grozza@sissa.it} }
%
%
\maketitle
\abstract*{In this manuscript,  starting with a Domain Decomposition (DD) algorithm on non-overlapping domains, we aim at the comparison of couplings of different discretisation models, such as Finite Element (FEM) and Reduced Order (ROM) models for separate subcomponents. 
In particular, we consider an optimisation-based DD model on two non-overlapping subdomains where the coupling on the common interface is performed by introducing a control variable representing a normal flux. Gradient-based optimisation algorithms are used to construct an iterative procedure to fully decouple the subdomain state solutions as well as to locally generate ROMs on each subdomain. Then, we consider FEM or ROM discretisation models for each of the DD problem components, namely, the triplet state1-state2-control.
We perform numerical tests on the backward-facing step Navier-Stokes problem to investigate the efficacy of the presented couplings in terms of optimisation iterations, optimal functional values and relative errors. }

\abstract{
Using Domain Decomposition (DD) algorithm on non--overlapping domains, we compare couplings of different discretisation models, such as Finite Element (FEM) and Reduced Order (ROM) models for separate subcomponents.
In particular, we consider an optimisation--based DD model where the coupling on the interface is performed using a control variable representing the normal flux. 
We use iterative gradient-based optimisation algorithms to decouple the subdomain state solutions as well as to locally generate ROMs on each subdomain. 
Then, we consider FEM or ROM discretisation models for each of the DD problem components, namely, the triplet state1--state2--control.
On the backward--facing step Navier-Stokes (NS) problem, we investigate the efficacy of the presented couplings in terms of optimisation iterations, optimal functional values and relative errors. }
    
\section{Introduction}
\label{sec:intro}

With the increased interest in complex problems, such as multiphysics and multiscale models, as well as real--time computations, there is a strong need for domain--decomposition (DD) segregated solvers and reduced--order models (ROMs). 
Segregated models decouple the subcomponents of the problems to use already existing state--of--the--art numerical codes in each component. 
These methods are extremely important for multi--physics problems when efficient subcomponent numerical codes are already available, or when we do not have direct access to the numerical algorithms for some parts of the systems. ROMs are very useful when dealing with real--time simulations or multi--query tasks. They split the computational effort into two stages: the offline stage, which contains the most expensive part of the computations, and the online stage, which performs fast computational queries using structures that were pre--computed in the offline stage.

The DD--ROM combination is very promising when dealing with different discretisation techniques on different subdomains. This small contribution inspired by~\cite{DECASTRO2023116398,  tezaur2023} aims at expanding the optimisation--based DD--ROM methods investigated in~\cite{prusak_stationary, prusak_nonstationary} to the hybrid numerical models, where separate subcomponents of the DD problem can be approximated by either a full order model based on Finite Element method (FEM) or a reduced order model. 
We refer to~\cite[Section 1.1]{DECASTRO2023116398} for an exhaustive overview of the literature on projection--based ROM--ROM and FEM--ROM couplings. In those works, the DD is primarily used as a tool for enhancing the efficiency of ROMs in complex scenarios and mostly Schwarz iterative methods are employed. 
Here, instead, we investigate a different optimisation--based approach that is particularly important in multi--physics contexts, where the stability of other DD methods might be compromised due to particular physical phenomena~\cite{prusak_thesis}.
Moreover, the optimisation--based coupling might be an extremely powerful tool when incorporating emerging very efficient non--intrusive data--driven models while making sure that surrogate DD models are physically consistent.


The rest of the paper is constructed as follows. 
In Section~\ref{sec:problem_formulation}, starting with a monolithic formulation of the discretised NS equations we describe an optimisation--based discrete DD model and we derive the optimal control problem.  
In Section~\ref{sec:rom_fem}, we discuss the ROM based on Proper Orthogonal Decomposition (POD) and the different FEM and/or ROM coupling techniques. 
Finally, in Section~\ref{sec:num_results} we provide numerical tests on the backward--facing step NS problem and draw some conclusions. 

\section{Problem formulation}
\label{sec:problem_formulation}
Consider the time-dependent incompressible NS equations, i.e.,
$$\partial_t u + u \cdot \nabla u - \nabla \cdot (\nu \nabla u) + \nabla p =0$$
subject to $\nabla \cdot u =0$. We start describing its monolithic formulation, then we introduce a discretised optimisation--based DD problem employing the implicit Euler time--stepping scheme and a FEM space discretisation. The resulting optimal control problem is set up at each time step, aiming at minimising the distance between the subdomain velocity fields by finding an optimal normal flux at the interface. 
\subsection{Monolithic formulation}
\label{monolithic_formulation}

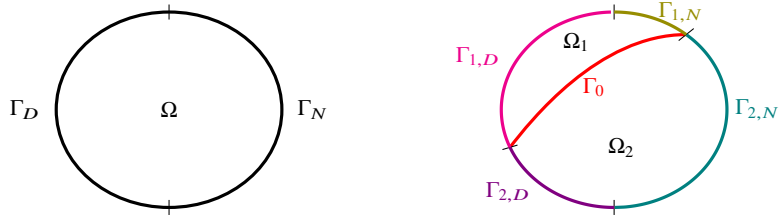
\begin{figure}\centering
\begin{subfigure}{0.49\textwidth}
       \centering
    \begin{tikzpicture}
 \draw[very thick] (0,0) ellipse (1.5cm and 1.3cm);
\node at (0,0) {$\Omega$};
\draw (0,1.2) -- (0,1.4);
\draw (0,-1.2) -- (0,-1.4);
\node[anchor=east] at (-1.6,0) {$\Gamma_{D}$};
\node[anchor=west] at (1.6,0) {$\Gamma_{N}$};
\end{tikzpicture}
\end{subfigure}\hfill
\begin{subfigure}{0.49\textwidth}
\centering
    \begin{tikzpicture}
\draw[red, very thick]  (0.9585,1) parabola (-1.3846,-0.5);
\node[red, very thick, anchor=east] at (-0., 0.3) {\textbf{$\Gamma_0$}};

\node at (-0.5,0.9) {$\Omega_1$};
\node at (0.1,-0.5) {$\Omega_2$};
\draw (0,1.2) -- (0,1.4);
\draw (0,-1.2) -- (0,-1.4);
\draw (1.0585,1.08) -- (0.8585,0.92);
\draw (-1.4846,-0.53) -- (-1.2846,-0.47);

\node[anchor=east, color=magenta] at (-1.4,0.7) {$\Gamma_{1, D}$};
\node[anchor=west, color=teal] at (1.5,0) {$\Gamma_{2, N}$};
\node[anchor=west, color=olive] at (0.48,1.3) {$\Gamma_{1, N}$};
\node[anchor=east, color=violet] at (-1,-1.1) {$\Gamma_{2, D}$};

\draw[domain=0:0.9585, smooth, very thick, variable=\x, samples=51, olive] plot ({\x}, {1.3*sqrt(1.0-\x*\x/1.5/1.5)});

\draw[domain=-1.3:1, smooth, very thick, variable=\y, samples=501, teal] plot ({sqrt(1.5*1.5*(1-\y*\y/1.3/1.3))}, {\y});

\draw[domain=-1.3846:0, smooth, very thick, variable=\x, samples=51, violet] plot ({\x}, {-1.3*sqrt(1.0-\x*\x/1.5/1.5)});

\draw[domain=-0.5:1.3, smooth, very thick, variable=\y, samples=1001, magenta] plot ( {-sqrt(1.5*1.5*(1-\y*\y/1.3/1.3))},{\y});

\end{tikzpicture}
    \end{subfigure}
    \caption{Domains and boundaries}\label{fig:dd_domain}
\end{figure}

Let $\Omega$ be a physical domain of interest: we assume $\Omega$ to be an open subset of $\mathbb{R}^2$ and $\Gamma$ to be the boundary of $\Omega$; see Fig.~\ref{fig:dd_domain} (on the left). We also consider a finite time interval $[0, T]$ with $T >0$.
Let $f: \Omega \times [0,T] \rightarrow \mathbb{R}^2$ be the forcing term, $\nu$ the kinematic viscosity, $u_{D}$ a given Dirichlet datum to be imposed on $\Gamma_D\subset \Gamma$ and $u_0$ a given initial condition.

Following~\cite{prusak_stationary, prusak_nonstationary}, we can define usual Lagrangian Finite Element (FE) spaces  on a triangulation $\Omega_h$ of $\Omega$ as follows: 
\begin{itemize}
    \item $V_{h} \subset \left[H^1(\Omega)\right]^2, \quad || \cdot||_{V_{h}} = || \cdot||_{\left[H^1(\Omega)\right]^2},  $
    \item $V_{0, h} \subset \left\{ v \in \left[H^1(\Omega)\right]^2: \left. v \right|_{\Gamma_{D} = 0} \right\}, \quad || \cdot||_{V_{0,h}} = || \cdot||_{\left[H^1(\Omega)\right]^2},  $
    \item $Q_{h} \subset L^2(\Omega), \quad || \cdot||_{Q_{h}} = || \cdot||_{L^2(\Omega)} $
\end{itemize}
and a time discretisation, through the time step $\Delta t >0$, and we assume the following time interval partition: $0 = t_0 < t_1 < .... < t_M =T$, where $t_n = n\Delta t$ for $n = 0, ..., M$.
 
The discretised problem with FEM and implicit Euler reads as follows at each time step $t_n$: find the velocity field $u_h^{n}: \Omega \rightarrow \mathbb{R}^2 $ and the pressure $p_h^{n}: \Omega\rightarrow \mathbb{R}$ s.t.
\begin{subequations}   \label{eq:mono}
\begin{align}
       \label{eq:state_fem_mono1}   \frac{ m ( u_{h}^n - u_{h}^{n-1}, v_{h})}{\Delta t}&+ a(u_{h}^n, v_{h}) + c(u_{h}^n, u_{h}^n, v_{h})     
  + b(v_{h}, p_{h}^n)    =   (f^n, v_{h})_{\Omega} 
     &\forall v_{h} \in V_{0,h},      \\
  \label{eq:state_fem_mono2} & b(u_{h}^n, q_{h})   = 0 \quad \quad  \forall q_{h} \in Q_{h},  \\ \label{eq:state_fem_mono3}
 & \quad \quad u^n    =   u_{D, h}^n  \quad \quad  \text{on}  \ \Gamma_{ D},
\end{align}
\end{subequations}
where 
\begin{itemize}
    \item  $m: V_{h} \times V_{0,h} \rightarrow \mathbb{R}, \quad m(u_{h},v_{h}) :=  (u_{h},  v_{h})_{\Omega}$,
    \item  $a: V_{h} \times V_{0,h} \rightarrow \mathbb{R}, \quad a(u_{h},v_{h}) := \nu (\nabla u_{h}, \nabla v_{h})_{\Omega} $,
     \item  $b: V_{h} \times Q_{h} \rightarrow \mathbb{R}, \quad b(v_{h},q_{h}) := - (\text{div} v_{h}, q_{h})_{\Omega}$,
     \item  $c: V_{h} \times V_{h} \times V_{0,h} \rightarrow \mathbb{R}, \quad c(u_{h},w_{h}, v_{h}) := \left( (u_{h} \cdot \nabla)w_{h}, v_{h}\right)_{\Omega}   $.
\end{itemize}
Moreover, since NS equations have a saddle--point structure, we require the pairs of spaces $V_{h}-Q_{h}$ and $V_{0,h}-Q_{h}$ to be inf--sup stable and this is achieved by using, for example, the Taylor--Hood $\mathbb P_2-\mathbb P_1$ FE spaces.



\label{DD}

\subsection{Discrete Domain Decomposition formulation}\label{sec:dd_fem}
As mentioned in the introduction, we resort to the optimisation--based approach for DD as described in~\cite{prusak_stationary, prusak_nonstationary}. Let $\Omega_i, \ i=1,2$, be two open subsets of $\Omega$, such that  $\overline{\Omega} = \overline{\Omega_1 \cup \Omega_2}$, \  $\Omega_1 \cap \Omega_2 = \emptyset$. Denote $\Gamma_i := \partial \Omega_i \cap \Gamma, \ i=1,2,$ and $\Gamma_0 := \overline{\Omega_1} \cap \overline{\Omega_2}$. In the same way, we define the corresponding boundary subsets $\Gamma_{i,D}$ and $\Gamma_{i, N}$, $i=1,2$, see Fig.~\ref{fig:dd_domain} (on the right). 


Next, following~\cite{prusak_stationary, prusak_nonstationary}, we assume to have at hand two well--defined triangulations $\mathcal T_1$ and $\mathcal T_2$ over the polygonal domains $\Omega_1$ and $\Omega_2$ respectively, and a one--dimensional discretisation $\mathcal T_0$ of the interface $\Gamma_0$.
We assume the meshes $\mathcal T_1$, $\mathcal T_2$ and $\mathcal T_0$ to share the same degrees of freedom (Dofs) on the interface $\Gamma_0$.

We can then define the restriction of the FE spaces $V_h,\,V_{h,0},\,Q_h$ onto the subdomain $\Omega_{i,h}$ as $V_{i,h},\,V_{0,i,h},\,Q_{i,h}$ and we define a FE space for the interface as
$$X_{h} \subset \left[L^2(\Gamma_0)\right]^2, \quad || \cdot||_{X_h} = || \cdot||_{\left[L^2(\Gamma_0)\right]^2} .$$
Similarly, we define the restriction of the bilinear and trilinear forms onto the restricted FE spaces as $m_i,\,a_i,\,b_i,\,c_i$. 

Also in the subcomponents, we choose the $\mathbb P_2-\mathbb P_1$ FE spaces to be inf-sup stable. For $X_h$, our choice is to use the $\mathbb P_2$ FE space implying that $X_h$ shares the DoFs with the spaces $V_{i,h}, i=1,2$ on the curve $\Gamma_0$.

The discretised optimisation--based DD formulation of the problem~\eqref{eq:mono} then reads as follows: for $n \geq 1$  
minimise over $g_h \in X_h$ the functional
  \begin{equation}
  \label{eq:functional_fem}  \mathcal J(u_{1,h}^n, u_{2,h}^n; g_h) = \frac{1}{2} \int_{\Gamma_0} \left| u_{1,h}^n - u_{2,h}^n\right|^2 d\Gamma \end{equation}
with each of the $u_{i,h}^n$ subject to the variational problem depending on $g_h$: \\
for $i=1,2$  find  $u_{i,h} \in V_{i,h}$  and $p_{i,h} \in Q_{i,h}$  satisfying 
\begin{subequations}\label{eq:state_fem}
\begin{align}
\begin{split}
       \label{eq:state_fem1}   \frac{ m_i ( u_{i,h}^n - u_{i,h}^{n-1}, v_{i,h})}{\Delta t}&+ a_i(u_{i,h}^n, v_{i,h}) + c_i(u_{i,h}^n, u_{i,h}^n, v_i)     \\ 
  + b_i(v_{i,h}, p_{i,h}^n)   & =   (f_i^n, v_{i,h})_{\Omega_i}
     + \left( (-1)^{i+1} g_h, v_{i,h} \right)_{\Gamma_0} 
\end{split}
     &\forall v_i \in V_{i,0,h},      \\
  \label{eq:state_fem2} & b_i(u_{i,h}^n, q_{i,h})   = 0, \quad \quad  &\forall q _{i,h} \in Q_{i,h}  \\ \label{eq:state_fem3}
 & \quad \quad u_i^n    =   u_{i, D, h}^n  \quad \quad  &\text{on}  \ \Gamma_{i, D},
\end{align}
\end{subequations}
where $u_{i,D,h}^n$ is the Galerkin projection of $u_{D}$ onto the trace--space $V_{i,h}|_{\Gamma_{i,D}}$. The control function $g_h$ entered the definition of the functional~\eqref{eq:functional_fem} implicitly through the solutions $u_{i,h}^n, i=1,2$ to~\eqref{eq:state_fem}. This is due the fact that, as proved in~\cite{prusak_nonstationary}, the optimal control problem~\eqref{eq:functional_fem}--~\eqref{eq:state_fem} is well--posed without the need of regularisation. 
\subsection{Optimality system and the gradient of the objective functional}
We now aim at providing the ingredients to set up a gradient--based iterative optimisation algorithm of the DD minimisation problem~\eqref{eq:functional_fem}--~\eqref{eq:state_fem}. 
In order to deal with variational constraint optimal control problem, we use the Lagrangian functional and sensitivity derivatives approaches; we refer to~\cite{prusak_stationary, prusak_nonstationary, Gunzburger_book, Hinze2009} for more details. 

The optimality system arising for the optimal control is defined as follows:
\begin{itemize}
    \item state problem~\eqref{eq:state_fem},
    \item adjoint problem: for $i=1,2$, find $\xi_{i,h} \in V_{i,0, h}$ and $\lambda_{i,h} \in Q_{i,h}$ satisfying
\begin{subequations}\label{eq:adjoint_fem}
\begin{align}
    \label{eq:adjoint1_fem} & &\frac{m_i(\eta_{i,h}, \xi_{i,h})}{\Delta t } + a_i(\eta_{i,h}, \xi_{i,h} ) +c_i (  \eta_{i,h} , u_{i,h}^n, \xi_i )  + c_i( u_{i,h}^n , \eta_{i,h}, \xi_{i,h} ) \\ \nonumber
& &+  b_i ( \eta_{i,h}, \lambda_{i,h}) =  ((-1)^{i+1}\eta_{i,h}, u_{1,h}^n - u_{2,h}^n)_{\Gamma_0}
\quad \forall \eta_{i,h} \in V_{i,0, h},
\\
 \label{eq:adjoint2_fem}  & &  b_i (\xi_{i,h}, \chi_{i,h})  =0  \quad \quad \quad   \forall \chi_{i,h}\in Q_{i,h}. 
\end{align}
\end{subequations} 
\item optimality condition:
\begin{equation}
\label{eq:optimality_condition_fem}    (r_h, \xi_{1.h} - \xi_{2,h})_{\Gamma_0} = 0  \quad \quad \forall r_h \in X_h.
\end{equation}
\end{itemize}

Resorting to the sensitivity derivatives approach \cite{prusak_stationary, prusak_nonstationary} leads to the following gradient representation of the objective functional~\eqref{eq:functional_fem} for a given $g_h \in X_h$:
\begin{equation}
  \label{eq:gradient_fem}  \frac{d\mathcal{J}}{dg}(u_{1,h}^n, u_{2,h}^n; g_h) =  \left.\xi_{1,h}\right|_{\Gamma_0} - \left.\xi_{2,h}\right|_{\Gamma_0},
\end{equation}
where $u_{i,h}^n, i=1,2$ are the solutions to the state equations~\eqref{eq:state_fem} and $\xi_{i,h}, i=1,2$ are the solutions to the adjoint equations~\eqref{eq:adjoint_fem}. Note that condition~\eqref{eq:optimality_condition_fem} ensures that the solutions to the coupled optimality system~\eqref{eq:state_fem}--~\eqref{eq:optimality_condition_fem} are the stationary points of the functional~\eqref{eq:functional_fem}.

\section{Reduced Order Model and FEM--ROM couplings}\label{sec:rom_fem}

As highlighted in Section~\ref{sec:intro}, reduced--order methods are efficient tools for significant reduction of the computational time for parameter--dependent PDEs.
This section deals with the reduced--order model for the problem obtained in the previous section, where the state equations, namely NS equations, are assumed to be dependent on a set of physical parameters. 
We study then different coupling options choosing for each subcomponent of the DD problem either the FEM or the ROM.


In this section, we will list all the necessary components to set--up a reduced order model. All the details can be found in~\cite{prusak_stationary, prusak_nonstationary}. 

Our goal is to generate linear low--dimensional subspaces of the FE spaces presented in Section~\ref{sec:dd_fem}. We rely on the POD compression technique; see, for instance,~\cite{prusak_stationary, prusak_nonstationary, Rozza_book}. The POD is based on the sampling of the parameter space $\mathcal P$ with a discrete set $\mathcal P_M$ and storing the snapshots associated with each parameter $\mu \in \mathcal P_M$ and each time instance. 
Since we aim at constructing linear spaces, we need to introduce parameter--dependent lifting functions $l_{i,h}^n(\mu) \in V_{i,h},$ for $ \mu \in \mathcal P_M$, such that $l_{i,h}^n(\mu) = u_{i,D,h}^n(\mu)$ on $\Gamma_{i,D}$ 
and define homogenised snapshots $u_{i,0,h}^n(\mu):=u_{i,h}^n(\mu) - l_{i,h}^n(\mu)$ which implies that $u_{i,0,h}^n(\mu)$ belongs to $V_{i,0,h}$. The reduced spaces $V_N^{u_1}, V_N^{p_1}, V_N^{u_2}, V_N^{p_2}$ and $V_N^{g}$ are then built as it is described in~\cite{prusak_stationary, prusak_nonstationary} together with the velocity supremiser technique~\cite{ballarin_supremiser}  in order to guarantee the inf--sup stability of the velocity--pressure reduced spaces.

Having at our disposal the reduced spaces, we perform the Galerkin projection of the state problem~\eqref{eq:state_fem}: for a given parameter $\mu \in \mathcal P$ and $g_N(\mu) \in V_N^g$, find {$u_{i,N}^n(\mu) = u_{i,0,N}^n(\mu) + l_{i,N}^n(\mu)$} with $u_{i,0,N}^n \in V_N^{u_i}$ and $p_{i,N}{(\mu)} \in V_N^{p_i} $ satisfying
\begin{subequations}\label{eq:state_rom}
\begin{align}
     \nonumber  \label{eq:state_rom1}   \frac{ m_i ( u_{i,N}^n - u_{i,N}^{n-1}, v_{i,N})}{\Delta t}&+ a_i(u_{i,N}^n, v_{i,N}) + c_i(u_{i,N}^n, u_{i,N}^n, v_{i,N})      \\ 
  + b_i(v_{i,N}, p_{i,N}^n)   & =   (f_i^n, v_{i,N})_{\Omega_i}
     + \left( (-1)^{i+1} g_N, v_{i,N} \right)_{\Gamma_0}
     &\forall v_{i,N} \in V_N^{u_1},      \\
  \label{eq:state_rom2} & b_i(u_{i,N}^n, q_{i,N})   = 0  &\forall q _{i,N} \in V_N^{p_i},
  \end{align}
\end{subequations}
where $l_{i,N}^n$ is the Galerkin projection of the lifting function $l_{i,h}^n$ to the finite dimensional space $V_N^{u_i}$ and $i=1,2$. In the previous equations, we omitted the dependence on the parameter $\mu$ in all forms (in particular, we will use the parametrised viscosity $\nu(\mu)$ in $a_i$ and the Dirichlet data $u^n_{i,D,h}(\mu)$) and all variables, for a clearer exposition.

Similarly, we can write the reduced counterpart of the adjoint equations~\eqref{eq:adjoint_fem}: for a given parameter $\mu \in \mathcal P$ and $u_{i,N}^N{(\mu)} \in V_N^{u_i} + \{l_{i,N}^n{(\mu)}\}$, find $\xi_{i,N} {(\mu)} \in V_N^{u_i}$ and $\lambda_{i,N}{(\mu)} \in V_N^{p_i} $ satisfying for $i=1,2$
\begin{subequations}\label{eq:adjoint_rom}
\begin{align}
\label{eq:adjoint1_rom} \frac{1}{\Delta t }m_i(\eta_{i,N}, \xi_{i,N}) &+ a_i(\eta_{i,N}, \xi_{i,N} ) +c_i \left(  \eta_{i,N} , u_{i,N}^n, \xi_i \right)  + c_i\left( u_{i,N}^n , \eta_{i,N}, \xi_{i,N} \right) \\ 
\nonumber& +  b_i ( \eta_{i,N}, \lambda_{i,h}) =  ((-1)^{i+1}\eta_{i,N}, u_{1,N}^n - u_{2,N}^n)_{\Gamma_0}, \ \ \forall \eta_{i,N} \in V_{i,N}^{u_i},
\\
 \label{eq:adjoint2_rom}  & \quad b_i (\xi_{i,N}, \chi_{i,N})  =0  ,  \quad \quad \quad \quad \forall \chi_{i,N}\in V_{i,N}^{p_i}. &
\end{align}
\end{subequations}

\subsection{FEM--ROM couplings: exploring different options}\label{sec:rom_fem_couplings}
In this section, we will provide a setting where different types of models, i.e. FEM or ROM,  can be chosen for each of the components of the DD problem, namely the triplet state1--state2--control. In particular, we investigate four different scenarios --- FEM--FEM--FEM, FEM--ROM--FEM, FEM--ROM--ROM and ROM--ROM--ROM. Each of these choices is characterised by a different optimisation problem. 

To proceed we need to define the following parameter--dependent projection and lifting operators from the reduced spaces onto the FE spaces:
\begin{itemize}
    \item $\Pi_i^n(\mu): V_{i,h} \rightarrow V_N^{u_i} + \{l_{i,N}^n(\mu)\} $, \quad $(\Pi_i^n(\mu))^T: V_N^{u_i}+\lbrace l_{i,N}^n(\mu)\rbrace  \rightarrow V_{i,h}$,
    \item $\Pi_{i,0}:  V_{i,0,h} \rightarrow V_N^{u_i}  $,\quad $\Pi_{i,0}^T: V_N^{u_i}  \rightarrow V_{i,0,h}$,
    \item $\Pi_X: X_{h} \rightarrow V_N^{g}$,\quad $\Pi_X^T: V_N^{g} \rightarrow X_{h} $.
\end{itemize}

In order to keep the exposition clear, we consider a generalised version of the objective functional \eqref{eq:functional_fem}
\begin{equation*}
   \mathcal J(u_1^n, u_2^n; g) = \frac{1}{2} \int_{\Gamma_0} \left| u_1^n - u_2^n\right|^2 d\Gamma 
\end{equation*}
with gradient
\begin{equation*}
     \frac{d\mathcal{J}}{dg}(u_1^n, u_2^n; g) =  \alpha_1 - \alpha_2,
\end{equation*}
where the state variables $u_1^n$, $u_2^n$, and the gradient contributions $\alpha_1$, $\alpha_2$ and the control $g$ assume the values reported in Table~\ref{tab:different_hybrid_choices}, according to the different hybrid method of choice. 
The minimisation problem, in abstract form, reads as follows: minimise over $g$ in the proper function space the functional
$\mathcal J(u_1^n, u_2^n; g)$. The expression of the corresponding gradient for each hybrid model can be inferred by looking at the values of the adjoint variables in Table~\eqref{tab:different_hybrid_choices}. We recall that $u_1^n$ and $u_2^n$ depend on $g$ in $\mathcal J$ through the FEM or ROM problems defined in the previous sections.
\begin{table}
    \centering
    \begin{tabular}{c|c|c|c|c}
         & {FFF} & {FRF} & {FRR} & {RRR}\\
         \hline
        $u_1^n$ & $u_{1,h}^n$  & $u_{1,h}^n$ & $\Pi_X u_{1,h}^n$ & $u_{1,N}^n$ \\
        \hline
        $u_2^n$ & $u_{2,h}^n$  & $(\Pi_2^n(\mu))^Tu_{2,N}^n$  & $\Pi_X (\Pi_2^n(\mu))^T u_{2,N}^n$ & $u_{2,N}^n$\\
        \hline
        $g$ & $g_h\in X_h$ & $g_h\in X_h$ & $g_N\in V_N^g$ & $g_N\in V_N^g$\\
        \hline
        $\alpha_1$ & $\xi_{1,h}|_{\Gamma_0}$  & $\xi_{1,h}|_{\Gamma_0}$ & $\Pi_X \xi_{1,h}$ & $\Pi_X\Pi_{1,0}^T\xi_{1,N}$\\
        \hline
        $\alpha_2$ & $\xi_{2,h}|_{\Gamma_0}$ & $\Pi_{2,0}^T\xi_{2,N}|_{\Gamma_0}$ & $\Pi_X\Pi_{2,0}^T\xi_{2,N}$  & $\Pi_X\Pi_{2,0}^T\xi_{2,N}$\\
        \hline
    \end{tabular}
    \caption{Different state, adjoint and control variables to use according to the different hybrid method in the DD approach. The method names are denoted by F for FEM and R for ROM for the variables $u_1$, $u_2$ and $g$ in this order. The solutions of $u^n_{i,h}$, $u_{i,N}^n$, $\xi_{i,h}$ and $\xi_{i,N}$ are obtained, respectively, in \eqref{eq:state_fem}, \eqref{eq:state_rom}, \eqref{eq:adjoint_fem} and \eqref{eq:adjoint_rom} }
    \label{tab:different_hybrid_choices}
\end{table}

\section{Numerical results}\label{sec:num_results}

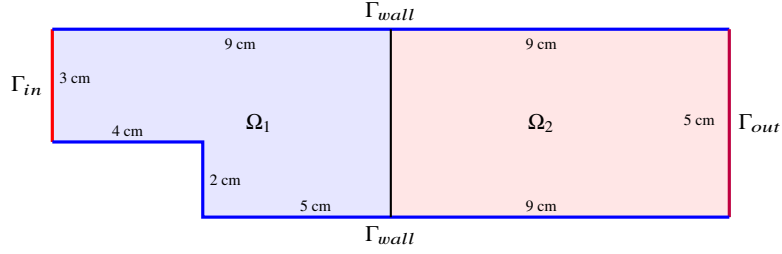
\begin{figure}
    \centering
\begin{tikzpicture}
    \fill[fill=blue!10] (0,1)  -- (0,2.5) -- (4.5,2.5) -- (4.5,0) -- (2,0)--(2,1)--(0,1); 
    \fill[fill=red!10] (4.5,0) -- (4.5,2.5)-- (9,2.5) -- (9, 1.25) -- (9,0) ;
    \draw[draw=red,very thick] (0,1) -- (0, 1.75) node [anchor=east]{$\Gamma_{in}$}  -- (0,2.5);
    \draw[draw=blue,very thick] (0,2.5) --  (4.5, 2.5) node [anchor=south] {$\Gamma_{wall}$}-- (9,2.5);
    \draw[draw=purple,very thick](9,2.5) -- (9, 1.25) node [anchor=west] {$\Gamma_{out}$}-- (9,0.);
    \draw[draw=blue,very thick] (9,0) -- (4.5, 0) node [anchor=north]{$\Gamma_{wall}$} -- (2,0)--(2,1)--(0,1);
    \draw[draw=black,thick](4.5,2.5) -- (4.5, 0);

    \node[black,scale=0.7] at (1,1.15) {4 cm};
    \node[black,scale=0.7] at (2.3,0.5) {2 cm};
    \node[black,scale=0.7] at (6.5,0.15) {9 cm};
    \node[black,scale=0.7] at (3.5,0.15) {5 cm};
    \node[black,scale=0.7] at (0.3,1.85) {3 cm};
    \node[black,scale=0.7] at (2.5,2.3) {9 cm};
    \node[black,scale=0.7] at (6.5,2.3) {9 cm};
    \node[black,scale=0.7] at (8.6,1.3) {5 cm};    
    \draw (2.75, 1.25) node {$\Omega_1$};
    \draw (6.5, 1.25) node {$\Omega_2$};
    
\end{tikzpicture}
\caption{Physical domain and domain decomposition for the backward--facing step problem\label{fig:domain_bfs}
}

         
         
\end{figure}

We consider the backward--facing step flow test case presented in~\cite{prusak_nonstationary}. Fig.~\ref{fig:domain_bfs} represents the physical domain of interest and the two--domain decomposition performed by dissecting the domain by a vertical segment.
In the offline phase, we consider two physical parameters: the viscosity $\nu \in [0.4,2]$ and the maximal magnitude $\bar U \in [0.5, 4.5]$ of the inlet profile $u_{in}(x,y) = \left(\bar U \cdot \frac{4}{9} (y-2)(5-y),0  \right)^T$ on $\Gamma_{in}=\lbrace (x,y): x= 0, y \in [2,5]\rbrace$. 
\begin{figure}[b]
    \centering
    \begin{subfigure}{0.49\textwidth}
        \includegraphics[width=\textwidth]{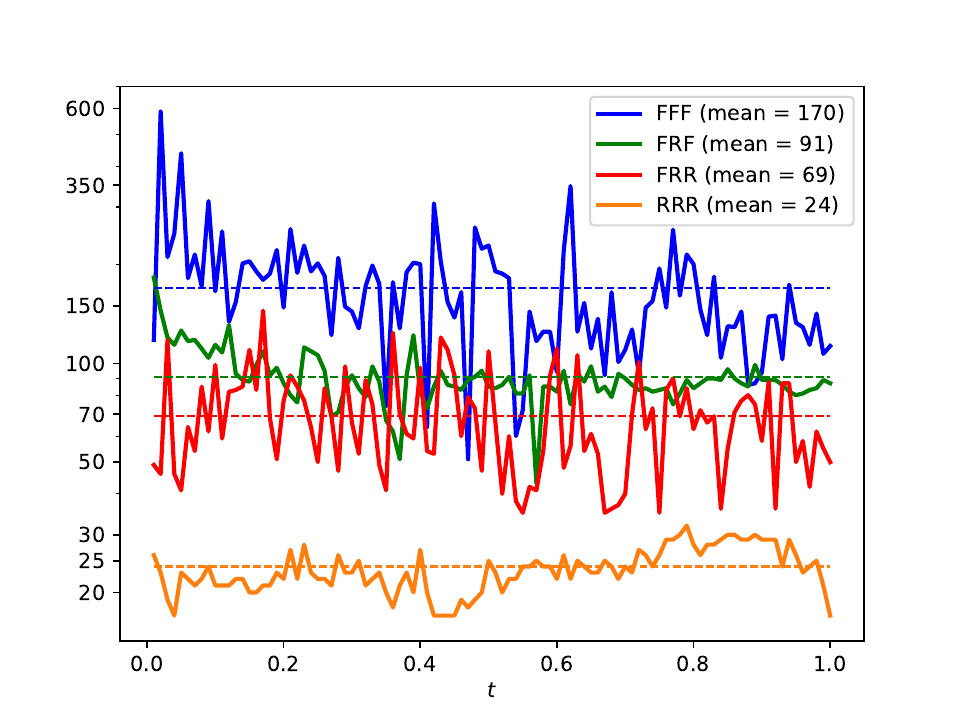}
        \caption{Optimisation iterations}\label{fig:iterations}
    \end{subfigure}
       \begin{subfigure}{0.49\textwidth}
        \includegraphics[width=\textwidth]{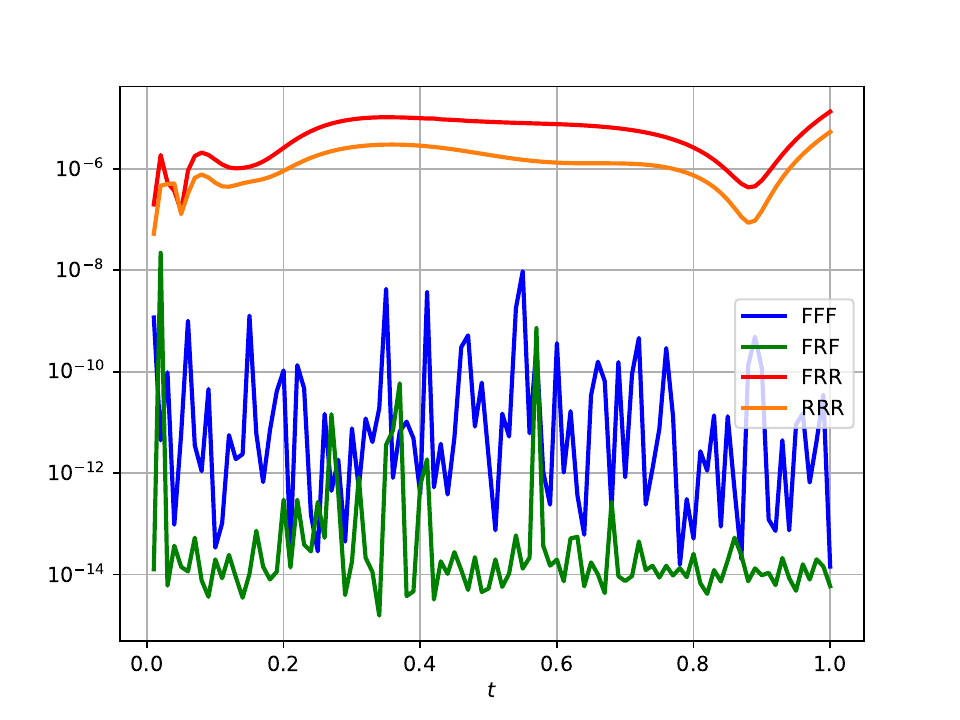}
        \caption{Functional values}\label{fig:func_values}
    \end{subfigure}
    \caption{Number of optimisation iterations (a) and the optimal functional values (b) of FFF, FRF, FRR and RRR solutions }\label{fig:its_func_values}
    \end{figure} 
    
In our test case, the FEM solutions are obtained on discrete state problems with a total of 27,890 DoFs, carrying out the minimisation in the interface space $X_h$ with 130 DoFs by the limited--memory Broyden--Fletcher--Goldfarb-–Shanno (L--BFGS--B) optimisation algorithm \cite{byrd1995limited} using the \texttt{scipy} library \cite{virtanen2020scipy}.
Snapshots are sampled from a training set of $64$ parameters randomly sampled from the 2--dimensional parameter space for $100$ time steps with $\Delta t = 0.01$ and $T=1$.
We perform a numerical analysis of the four couplings described in Section~\ref{sec:rom_fem_couplings} for a test parameter value $\left(\bar U, \nu\right) = \left(4.5, 0.4\right)$. For the couplings using the ROM model, we choose the following number of reduced basis modes: 30 for $u_1$, 12 for $u_2$, and 5 for each of $p_1$ and $p_2$ and the corresponding supremisers.
    \begin{figure}[b]
        \centering
        \begin{subfigure}{0.49\textwidth}
    \includegraphics[width=\textwidth]{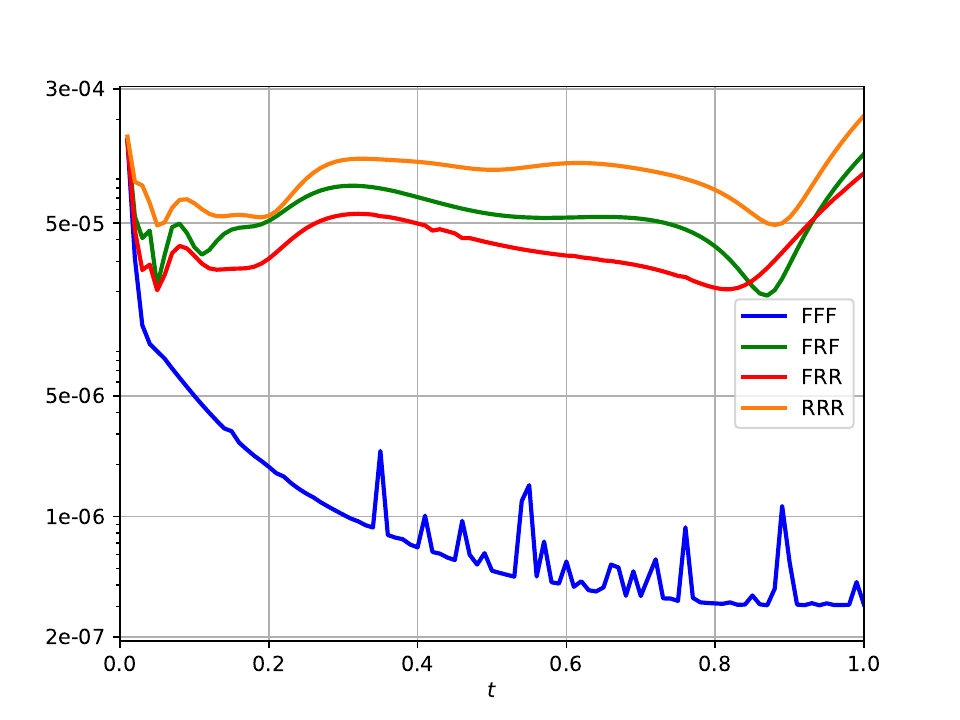}  
    \caption{Relative error for $u_1$}\label{fig:rel_u1}
    \end{subfigure}
    \begin{subfigure}{0.49\textwidth}
    \includegraphics[width=\textwidth]{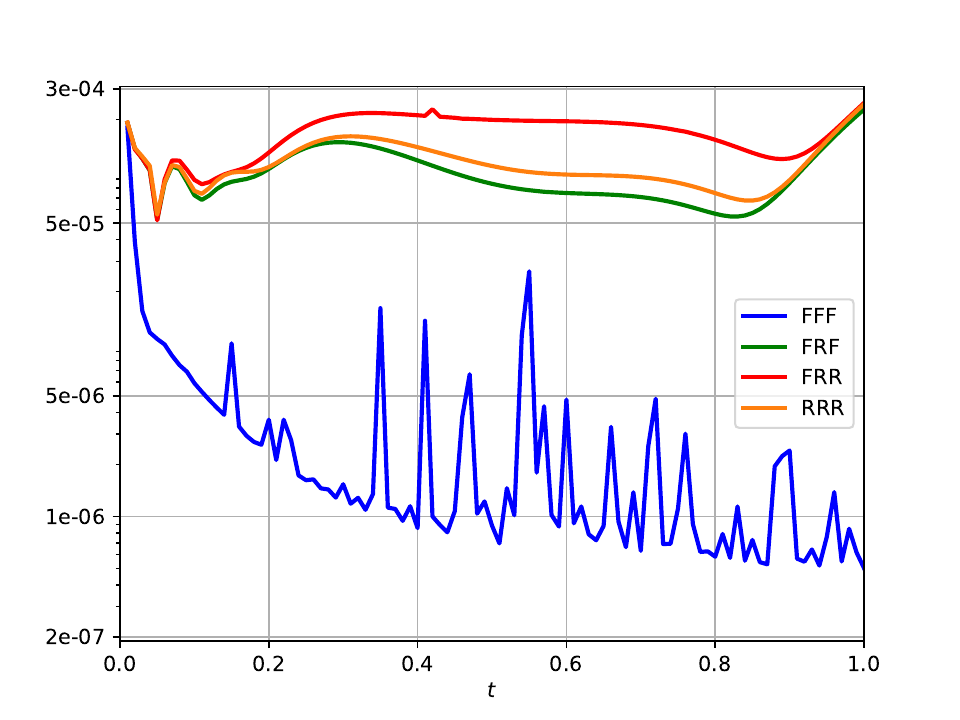}  
    \caption{Relative error for $u_2$}\label{fig:rel_u2}
    \end{subfigure}\\
        \begin{subfigure}{0.49\textwidth}
    \includegraphics[width=\textwidth]{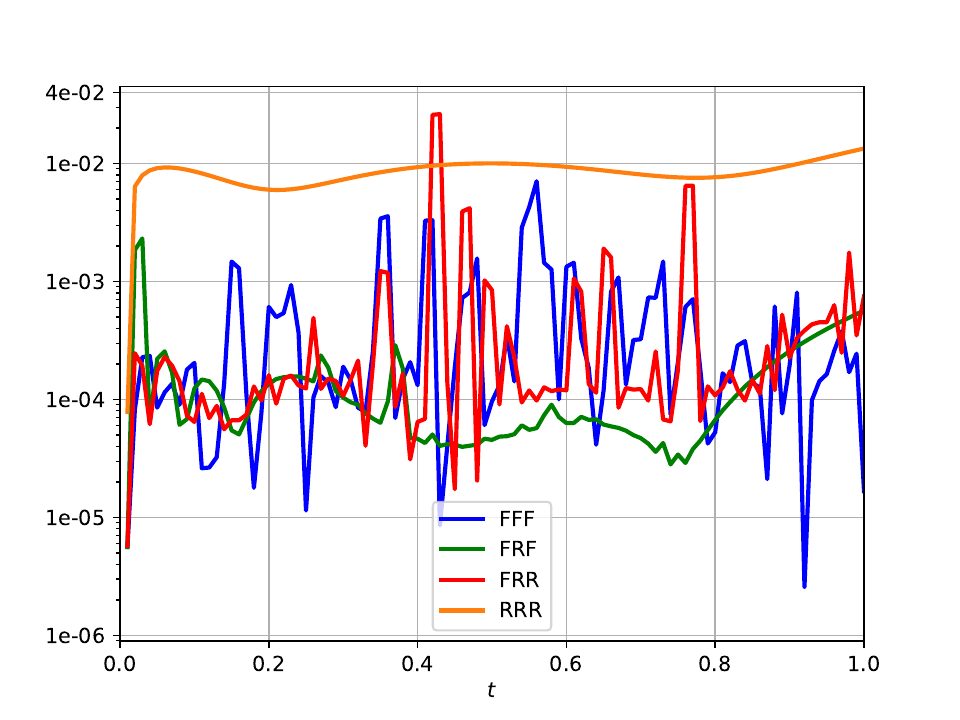}  
    \caption{Relative error for $p_1$}\label{fig:rel_p1}
    \end{subfigure}
    \begin{subfigure}{0.49\textwidth}
    \includegraphics[width=\textwidth]{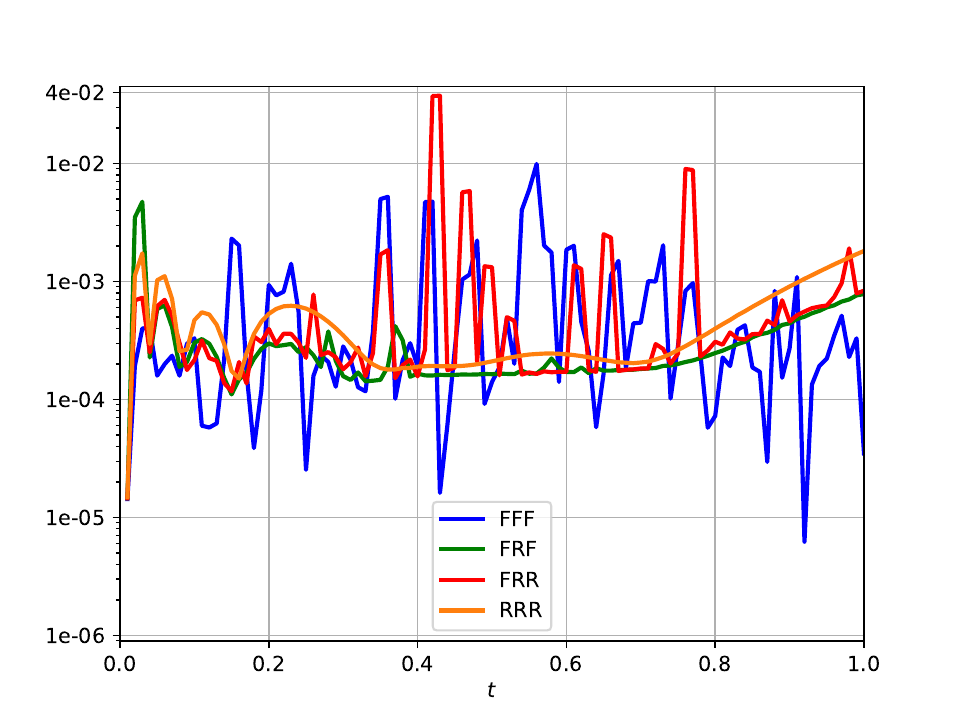}  
    \caption{Relative error for $p_2$}\label{fig:rel_p2}
    \end{subfigure}
    \caption{Relative errors of FFF, FRF, FRR, and RRR solutions w.r.t. the monolithic solution}
    \label{fig:rel_errors}
\end{figure}

Fig.~\ref{fig:iterations} shows the number of optimisation iterations over time for couplings FFF, FRF, FRR and RRR, and Fig.~\ref{fig:func_values}  shows the relative error of each state subcomponent with respect to the monolithic solution of the NS problem~\eqref{eq:mono}. 
It can be easily seen that the FFF coupling has the overall highest number of iterations while the RRR coupling has the lowest, and the FRF and FRR are in between the other two. This is because the optimisation in the case of ROM model for the control variable is carried out over a much smaller set of admissible solutions, which is constructed on preliminary physical information, i.e., FEM snapshots. 
On the other hand, as expected, this is balanced by the relative errors as it is shown in Fig.~\ref{fig:rel_errors}, where the relative error in the FFF coupling scenario is much smaller (at least for the velocity fields $u_1$ and $u_2$) than for other types of coupling. 
The irregular nature of the relative errors for the pressure fields $p_1$ and $p_2$, which are nevertheless sufficiently low, is most probably due to the fact that we use a ``black--box'' optimisation algorithm. 
Indeed, we do not have access to how the search direction for the control variable is chosen. 
We believe this can be improved by the use of different minimisation algorithms where we can have more control over the iterative procedure and impose some additional constraints also on the pressures.
{As it can be seen from the optimal functional accounting for the interface error in Fig.~\ref{fig:func_values} and from Fig.~\ref{fig:rel_u2}, the RRR coupling outperforms the FRR coupling. Indeed, the FRR optimisation process shows convergence issues and it stagnates, leading to an excessive number of functional and gradient evaluations for defining a search direction. Even though the FRR method is able to keep relative errors quite low, this issue, most probably due to the incapacity of capturing the high--dimensional FEM dynamics at the interface by a handful of the reduced interface space DoFs, makes it unreliable. Further comprehensive theoretical analysis of the resulting discrete optimisation problems is needed and it will be the subject of future work. }

{Overall, we can conclude that already choosing one state variable to be in the reduced space can alleviate the costs of the optimisation by a factor of 2. This might be helpful when in a subdomain one might need extra accuracy while in another a ROM is enough. In that situation we suggest to use the FRF instead of the FRR model for the above mentioned issues. } 

{Lastly, it is worth mentioning that for a relatively simple test case as the one considered in the paper, monolithic FEM and ROM solutions are more computationally efficient w.r.t. any DD alternative considered here, as the added cost of iterative procedure dominates the overall cost of monolithic solvers. However, such hybrid models are relevant in the case of complex problems (especially multiphysics problems) that require a larger number of DoFs and sometimes are extremely difficult (or even impossible) to confront in a feasible computational time.   }
\begin{acknowledgement}
M.~N. has been funded by the Austrian Science Fund (FWF) \href{https://doi.org/10.55776/ESP519}{10.55776/ESP519}. 
\end{acknowledgement}

\end{document}